\input amstex.tex
\documentstyle{amsppt}
\NoRunningHeads
\magnification=\magstep1
\baselineskip=12pt
\parskip=5pt
\parindent=18pt
\topskip=10pt
\leftskip=0pt
\pagewidth{32pc}
\pageheight{47pc}
\topmatter
\title The Simple Group of Order 168 and K3 Surfaces
\endtitle
\author Keiji Oguiso and De-Qi Zhang
\endauthor
\address
Department of Mathematical Sciences, University of Tokyo,
Komaba, Meguro, Tokyo Japan;
Department of Mathematics, National University of Singapore,
2 Science Drive 2, Singapore 117543
\endaddress
\email
oguiso$\@$ms.u-tokyo.ac.jp; matzdq$\@$math.nus.edu.sg
\endemail
\dedicatory
Dedicated to Herr Professor Doctor Hans Grauert on the occasion of his seventieth birthday
\enddedicatory
\abstract
The aim of this note is to characterize a K3 surface of Klein-Mukai type
in terms of its symmetry.
\endabstract
\endtopmatter
\document
\head Introduction \endhead 
The group $L_{2}(7)$ is by the
definition the projectivized special linear group $PSL(2, \bold
F_{7})$ and is generated by the three projective transformations
of $\bold P^{1}(\bold F_{7})$ of order $7$, $3$, $2$:
$$\alpha : x \mapsto x + 1; \beta : x \mapsto 2x; \gamma : x \mapsto -x^{-1},$$ where the coefficient $2$ in $\beta$ is a generator of the cyclic group
$(\bold F_{7}^{\times})^{2} (\simeq \mu_{3})$. (See for instance 
[CS, Chapter 10].) As well-known, this group is
of order
$168$ and is characterized as the second smallest non-commutative simple group.
\par
\vskip 4pt
One of
interesting connections
between $L_{2}(7)$ and complex algebraic geometry goes back to the result
of the great German mathematicians Hurwitz and Klein in G\"ottingen: 
$\vert L_{2}(7) \vert = 84(3-1)$ is the largest possible
order
of a group acting on a genus-three curve and the
so called Klein quartic curve
$$C_{168} =
\{x_{1}x_{2}^3 + x_{2}x_{3}^3 + x_{3}x_{1}^3 =  0\} \subset \bold P^{2}$$
is the unique genus-three curve admitting an $L_2(7)$-action.
The action of $L_{2}(7)$ on $C_{168}$ is
the projective transformation induced by
(one of two essentially the same) $3$-dimensional irreducible representation 
$V_{3}$ of $L_{2}(7)$ given by  
$$\alpha \mapsto  \pmatrix \zeta_{7} & 0 & 0 \\
               0 & \zeta_{7}^{2} &  0 \\
               0 &  0 & \zeta_{7}^{4}\\
               \endpmatrix ; \,\,\,\,
\beta \mapsto \pmatrix 0 & 0 & 1 \\
               1 & 0 &  0 \\
               0 & 1 & 0 \\
               \endpmatrix ; \,\,\,\,
\gamma \mapsto \frac{-1}{\sqrt{-7}} \pmatrix a & b & c \\
               b & c & a \\
               c & a & b \\
               \endpmatrix, $$
where $\zeta_7 = \text{exp} (2 \pi \sqrt{-1}/7)$, 
$a = \zeta_{7}^{2} - \zeta_{7}^{5}$, $b = \zeta_{7} - \zeta_{7}^{6}$,
$c = \zeta_{7}^{4} - \zeta_{7}^{3}$, and the branch of $\sqrt{-7}$ is chosen 
so that $\sqrt{-7} =  \zeta_{7} + \zeta_{7}^{2} + \zeta_{7}^{4} - \zeta_{7}^{3} 
- \zeta_{7}^{5} - \zeta_{7}^{6}$. The other $3$-dimensional irreducible 
representation is the composition of the representation $V_{3}$ with the outer
automorphism of $L_2(7)$ given by $\alpha \mapsto \alpha^{-1}$,
$\beta \mapsto \beta$ and $\gamma \mapsto \gamma$ (see [ATLAS] and
[Bu, Section 267]).
The Klein curve $C_{168}$ together with $L_{2}(7)$-action also appears in the 
McKay correspondece problem [Ma] and a classification of Calabi-Yau 
threefolds [Og].
\par 
\vskip 4pt
Our interest in this note is a relation between $L_{2}(7)$ and K3
surfaces. 
\par 
\vskip4pt 
Throughout this note, a {\it K3 surface} means a simply-connected 
smooth complex {\it algebraic} surface $X$ with a nowhere vanishing 
holomorphic 2-from $\omega_{X}$. We call an automorphism $g \in \text{Aut}(X)$ 
{\it symplectic} if $g^{*}\omega_{X} = \omega_{X}$. According to Mukai's 
classification [Mu1],
there are eleven maximum finite groups acting on K3 surfaces
symplectically, and among them,
there appear two simple groups: the group $L_{2}(7)$ and the alternating group
$A_{6}$ of degree $6$.
In the same paper, Mukai also gives a beautiful example of K3 surface $X_{168}$ with
$L_{2}(7)$-action, where
$$X_{168} = \{x_{0}^{4} +  x_{1}x_{2}^3 + x_{2}x_{3}^3
+ x_{3}x_{1}^3 =  0\} \subset \bold P^{3}.$$ This is the cyclic cover of
$\bold P^{2}$ of degree 4 branched along the Klein quartic curve $C_{168}$.
Here, the action of $L_{2}(7)$ on $X_{168}$ is
naturally induced by the action on $C_{168}$.
Note that this $X_{168}$ now admits a larger group action of
$L_{2}(7) \times \mu_{4}$,
where $\mu_{4}$ is the Galois group of the covering.
On the other hand, the smooth plane curve
$H_{168}$ of degree 6 defined by
$\{5x_{1}^{2}x_{2}^{2}x_{3}^{2} -
x_{1}^{5}x_{2} - x_{2}^{5}x_{3} - x_{3}^{5} x_1 = 0\}
\subset \bold P^{2}$ -- the zero locus of the Hessian of
the Klein quartic curve -- is also
invariant under the same $L_{2}(7)$-action on $\bold P^{2}$.
So, the K3 surface
$$X_{168}' =
\{y^{2} = 5x_{1}^{2}x_{2}^{2}x_{3}^{2} -
x_{1}^{5}x_{2} - x_{2}^{5}x_{3} - x_{3}^{5}x_{1}\}
\subset \bold P(1,1,1,3),$$
i.e.
the double cover of $\bold P^{2}$ branched along $H_{168}$,
also admits an $L_{2}(7)$-action.
However, it will turn out that these two K3 surfaces $X_{168}$ and $X_{168}'$
are
not isomorphic to each other (see {\bf Remark (2.12)}).
Therefore, K3 surfaces having
$L_{2}(7)$-action are no more unique and it is of interest to characterize
the Klein-curve-like K3 surface $X_{168}$
in a flavour similar to that of Hurwitz and Klein.
This is the aim of this short note.

\par \vskip 4pt
Throughout this note, we set $G := L_{2}(7)$. Our main observation is 
as follows:

\proclaim{Main Theorem} Let $X$ be a $K3$ surface. Assume that 
$G \subset \text{Aut}(X)$. Let $\tilde{G} $ be a finite subgroup of 
$\text{Aut}(X)$ such
that $G \subset \tilde{G}$. Then, 
\roster 
\item $\tilde{G}/G$ is a cyclic group of
order $1$, $2$, $3$, or $4$; and 
\item if  $\tilde{G}/G$ is of the maximum order $4$, then $(X, \tilde{G})$ is
isomorphic to
the Klein-Mukai pair $(X_{168}, L_{2}(7) \times \mu_{4})$. 
\endroster
\endproclaim
Here an isomorphism means an equivariant isomorphism with respect to group
actions.

\par \vskip 4pt
The main difference between genus-three curves and K3 surfaces is that there
are no canonical polarizations on K3 surfaces. In other words, we do not know
{\it a priori} which K3 surfaces are quartic K3 surfaces
or which polarizations are invariant under the group action.
Indeed, the determination of the invariant polarization for $\tilde{G}$
-- this will turn out to be of degree four if $\vert \tilde{G}/G \vert = 4$ 
({\bf Claim (2.10)}) --
is the most crucial part in this note.
Besides Mukai's pioneering work, we are much inspired by
a series of Kondo's work [Ko1,2]  on a lattice theoretic proof of Mukai's
classification and the determination of the K3 surface with the largest
finite group action as well as the action. Especially we will fully exploit
his brilliant idea of studying invariant lattices through an embedding of
their orthogonal complements into some Niemeier lattices.
This enables us to relate the problem with the Mathieu group $M_{24}$ and
the binary Golay code $\Cal C_{24}$ (Section one) and provides a very
powerful tool in calculating the discriminants of
the invariant lattices $H^{2}(X, \bold{Z})^{G}$ also in our setting. Combining
this with the additional
group action $\tilde{G}/G$, we shall determine the invariant polarization in
the maximum case $\vert \tilde{G}/G \vert = 4$. Once we find the invariant
polarization
in a lattice-theoretic way,
we can continue the proof by coming back to more algebro-geometric arguments.
One of the advantages of the algebro-geometric argument is perhaps
that we can then express the K3 surface and the group action
in a very concrete way as in the Theorem.

\par \vskip 4pt
\head Acknowledgement 
\endhead
Both authors would like to express their thanks to Professor M. L. Lang
for his help on Mathieu groups, to Professor Dolgachev for his critical reading of the first draft and pointing out the reference [Ed]
for the use in the proof of Claim 2.11, 
and to Professors S. Mukai and S. Kondo for their valuable comments. The main 
idea of this work was found during the 
first named author's stay at the National University of Singapore in August
2000.
He would like to express his thanks to the National
University of Singapore for financial support and to staff members there
for the warm hospitality.
\head 1. The Niemeier Lattices \endhead 
In this section, we recall some basic facts on the Niemeier 
lattices needed in our arguments.
Our main reference concerning Niemeier lattices and their relations
with Mathieu groups is [CS, Chapters 10, 11, 16, 18]. 
\definition{(1.1)} In this note, the even {\it negative} definite unimodular lattices of rank 24
are called Niemeier lattices. (We changed the sign from positive into
negative.)
There are exactly 24 isomorphism classes of the Niemeier lattices and
each isomorphism class is uniquely determined by its root lattice $N_{2}$, 
i.e. the sublattice generated by all the roots, the elements $x$ with
$x^{2} = -2$.
Except the so called Leech lattice which contains no roots,
the other 23 lattices are the over-lattices of their root lattices, 
which are:
$$\gather A_{1}^{\oplus 24}, \, A_{2}^{\oplus 12}, \, A_{3}^{\oplus 8}, \,
A_{4}^{\oplus 6}, \, A_{6}^{\oplus 4}, \, A_{8}^{\oplus 3}, \,
A_{12}^{\oplus 2}, \, A_{24}, \\
D_{4}^{\oplus 6}, \, D_{6}^{\oplus 4}, \, D_{8}^{\oplus 3},
\, D_{12}^{\oplus 2}, \,
D_{24}, \, E_{6}^{\oplus 4}, \, E_{8}^{\oplus 3}, \, A_{5}^{\oplus 4}
\oplus D_{4},
A_{7}^{\oplus 2} \oplus D_{5}^{\oplus 2}, \\
A_{9}^{\oplus 2} \oplus D_{6}, \, A_{15} \oplus D_{9}, \, E_{8}
\oplus D_{16}, \,
E_{7}^{\oplus 2} \oplus D_{10}, \, E_{7} \oplus A_{17}, \,
E_{5} \oplus D_{7} \oplus A_{11}.
\endgather $$
We denote the Niemeier lattices $N$ whose root lattices are
$A_{1}^{\oplus 24}$, $A_{2}^{\oplus 12}$
and so on by $N(A_{1}^{\oplus 24})$ ,  $N(A_{2}^{\oplus 12})$ and so on. 
\enddefinition 
\definition{(1.2)} In what follows,
we regard the set of roots $R := \{r_{i} \vert 1 \leq i \leq 24\}$ 
corresponding to the vertices of the
Dynkin diagram, as the set of the simple roots of $N$.
Denote by $O(N)$ (resp. by $O(N_{2}))$ the group of isometries of $N$ 
(resp. of $N_{2}$) and by $W(N)$ the
Weyl group generated by the reflections given by the roots of $N$. 
Here $O(N) \subset O(N_{2})$ and $W(N)$ is a normal subgroup of both $O(N)$ and $O(N_{2})$. The invariant hyperplanes of
the reflections divide $N \otimes \bold R$ into (finitely many) chambers.
Then, each chamber is a fundamental domain of the action of $W(N)$ and
the quotient group $S(N) := O(N)/W(N)$ is identified with a subgroup of 
symmetry of
the distinguished chamber ${\Cal C} :=
\{x \in N \otimes \bold R \vert (x, r) > 0 r \in R \}$ and also a subgroup of 
a larger group $S_{24} = \text{Aut}_{\text{set}}(R)$. 
\enddefinition 
The groups $S(N)$ are very explicitly calculated in [CS, Chapters 18, 16].
(See also [Ko1].)  The following is a part of the results there:
\proclaim{Lemma (1.3) [CS, Chapters 18, 16]} Let $N$ be a non-Leech 
Niemeier lattice. Then,  
\roster 
\item 
$S(N) = M_{24}$ if
$N = N(A_1^{\oplus 24});$
\item $S(N) = C_2 \rtimes  (C_2^{\oplus 3} \rtimes  L_3(2))$ if
$N = N(A_3^{\oplus 8})$; and 
\item for other $N$, the order
$|S(N)|$ is not divisible by $7$. \qed
\endroster
\endproclaim 
Let us add a few remarks about the groups appearing in the Lemma above. 
The next (1.4) and (1.5) are concerned with the first case (1) and (1.6) is 
for the second case (2).
\definition{(1.4)}  
Observe that $A_{1}^{\oplus 24} \subset 
N(A_{1}^{\oplus 24}) \subset (A_{1}^{\oplus 24})^{*}$ and that 
$$(A_{1}^{\oplus 24})^{*}/A_{1}^{\oplus 24} = \oplus_{i=1}^{24} \bold{F}_{2} 
\overline{r}_{i} \simeq \bold{F}_{2}^{\oplus 24}.$$ 
Here $\overline{r}_{i} := r_{i}/2 \, \text{mod} \bold{Z}r_{i}$ is the standard 
basis of the i-th factor $(A_{1})^{*}/A_{1}$. We also identify 
$$S_{24} = \text{Aut}_{\text{set}}(R) = 
\text{Aut}_{\text{set}}(\{\overline{r}_{i}\}_{i = 1}^{24}).$$
The linear subspace of $(A_{1}^{\oplus 24})^{*}/A_{1}^{\oplus 24}$ 
$$\Cal C_{24} := N(A_{1}^{\oplus 24})/A_{1}^{\oplus 24} \simeq 
\bold{F}_{2}^{\oplus 12}$$ 
encodes the information which elements of $(A_{1}^{\oplus 24})^{*}$ lie in 
$N(A_{1}^{\oplus 24})$. Besides this role, this subspace $\Cal C_{24}$ 
carries the structure of the binary self-dual code of Type $II$ with minimal 
distance 8, called the (extended) binary Golay code. 
\par
Among many equivalent definitions, the Mathieu group $M_{24}$ of degree 
24 is defined to be the subgroup of 
$S_{24}$ preserving $\Cal C_{24}$, i.e.  
$$M_{24} := \{\sigma \in S_{24} \vert \sigma(\Cal C_{24}) = \Cal C_{24}\}.$$  
As well-known, $M_{24}$ is a simple group of order $24 \cdot 23 \cdot 22 \cdot 
21 \cdot 20 \cdot 16 \cdot 3$ and acts on the set 
$\{\overline{r}_{i}\}_{i = 1}^{24}$ as well as on $R$ quintuply transitively. 
\par  
Let $\Cal P(R)$ be the power set of $R$, i.e. the set consisting of the subsets 
of $R$. Then, $\Cal P(R)$ bijectively corresponds to the set 
$(A_{1}^{\oplus 24})^{*}/A_{1}^{\oplus 24}$ by: 
$$\iota : \Cal P(R) \ni A \mapsto \overline{r}_{A} := 
\frac{1}{2} \sum_{r_{j} \in A}r_{j} 
\, \text{mod} \, A_{1}^{\oplus 24} 
\in (A_{1}^{\oplus 24})^{*}/A_{1}^{\oplus 24}.$$ 
Set $\Cal E := \iota^{-1}(\Cal C_{24})$. 
Then $A \in \Cal E$ if and only if $\frac{1}{2} \sum_{r_{j} \in A}r_{j}$ is in 
$N(A_{1}^{\oplus 24})$. Moreover, it is known that $\emptyset, R \in \Cal E$ 
and that if $A \in \Cal E$ ($A \not= R$, $\emptyset$) then $\vert A \vert$ is 
either $8$, $12$, or $16$. We call $A \in \Cal E$ an Octad (resp. a Dodecad) 
if $\vert A \vert = 8$ (resp. $12$). Note that $B \in \Cal E$ with 
$\vert B \vert = 16$ is then the complement of an Octad (in $R$), i.e. 
$B$ is of the form $R - A$ for some Octad $A$. There are exactly 759 Octads. 
\enddefinition 

The following fact called the Steiner property $St(5, 8, 24)$ and its proof 
are both needed in the proof of our main result:
 
\proclaim{Fact (1.5) (the Steiner property)} For each 5-elemet subset 
$S$ of $R$, there exists exactly one Octad $O$ such that 
$S \subset O$. \endproclaim 
 
\demo{Proof} Since $M_{24}$ is quintuply transitive on $R$, there 
exists an Octad $O$ such that $S \subset O$. Let $O_{1}$ and $O_{2}$ 
be two Octads. Then, by the definition, their symmetric difference 
$(O_{1} - O_{2}) \cup (O_{2} - O_{1})$ is also an element of $\Cal E$. 
Thus, $\vert (O_{1} - O_{2}) 
\cup (O_{2} - O_{1}) \vert$ is either $0$, $8$, $12$, or $16$ and we have that 
$\vert O_{1} \cap O_{2} \vert$ is either $8$, $4$, $2$ or $0$. 
Therefore, if $S \subset O_{1}$ and $S \subset O_{2}$, then one has 
$\vert O_{1} \cap O_{2} \vert \geq 5$, whence 
$\vert O_{1} \cap O_{2} \vert = 8$. This means that $O_{1} = O_{2}$. \qed 
\enddemo           
\definition{(1.6)}
In the second case, we identify (non-canonically) the set of eight connected 
components of the Dynkin diagram
$A_{3}^{\oplus 8}$ with the three-dimensional linear space
$\bold F_{2}^{\oplus 3}$
over $\bold F_{2}$ by letting one connected component to be $0$. The group
$C_{2} \rtimes (C_{2}^{\oplus 3} \rtimes L_3(2))$ is the semi-direct product,
where $C_{2}$ interchanges the two edges of all the components,
$C_{2}^{\oplus 3}$ is the group of the parallel transformations of the affine
space
$\bold F_{2}^{\oplus 3}$ and $L_{3}(2) (\simeq L_{2}(7))$ is the linear
transformation group of $\bold F_{2}^{\oplus 3}$ which fixes
(point wise) the three simple roots in the identity component. \enddefinition 

\head 2. Proof of the main Theorem \endhead
In what follows, we set $L := H^{2}(X, \bold Z)$,
$L^{G} := \{x \in L \vert g^{*}x = x$
for all $g \in G\}$, and $L_{G} := (L^{G})^{\perp}$ in $L$.
This $L$ is the unique even unimodular lattice of index $(3, 19)$.
We also denote by $S_{X}$ the N\'eron-Severi lattice and by
$T_{X}$ ($:= S_{X}^{\perp}$ in $L$), the transcendental lattice.
Since $G$ is simple and non-commutative, we have $G = [G, G]$.
In particular, $G$ acts on $X$ symplectically. Therefore $L^{G}$ contains
both $T_{X}$
and the invariant ample classes under $G$, namely the pull back of ample
classes of $X/G$.
In addition, since $G$ is maximum among finite symplectic group actions [Mu1],
$G$ is normal in $\tilde{G}$ and the quotient group $\tilde{G}/ G$ acts
faithfully on
$H^{2,0}(X) = \bold C \omega_{X}$. In particular, $\tilde{G}/G$ is a cyclic
group of
order $I$ such that the Euler function $\varphi(I)$ divides rank $T_{X}$ [Ni].

\proclaim{Claim (2.1)}
\roster
\item $\text{\rm rank}\, L^{G} = 3$. In particular,
$\text{\rm rank} \, T_{X} = 2$ and
(up to scalar) there is exactly one $G$-invariant algebraic cycle class $H$.
Moreover, this class $H$ is ample and is also invariant under $\tilde{G}$.      \item $\vert \tilde{G}/G \vert$ is either $1$, $2$, $3$, $4$ or $6$.
\endroster
\endproclaim
\demo{Proof} The equality $\text{rank}\, L^{G} = 3$ is a special case of 
a general formula of Mukai. Here for the convenience of readers, we shall 
give a direct argument along [Mu1, Proposition 3.4]. Let us consider the 
natural representaion $\rho$ of $G$ on the cohomology ring of $X$ 
$$\tilde{L} := \oplus_{i= 0}^{4} H^{i}(X, \bold{Z}) = H^{0}(X, \bold{Z}) 
\oplus L 
\oplus H^{4}(X, \bold{Z}).$$ 
Then, by the representation theory of finite groups and by the 
Lefschetz $(1,1)$-Theorem, one has  
$$2 + \text{\rm rank}\, L^{G} = \text{\rm rank}\, \tilde{L}^{G} = 
\frac{1}{\vert G \vert} 
\sum_{g \in G} \text{tr}(\rho(g)) = \frac{1}{\vert G \vert} 
\sum_{g \in G} \chi_{\text{top}}(X^{g}).$$ 
Here, the terms in the last sum are calculated by Nikulin [Ni] as follows:  
$$\chi_{\text{top}}(X^{g}) = 24, 8, 6, 4, 4, 2, 3, 2$$ if 
$$\text{ord}(g) = 1, 2, 3, 4, 5, 6, 7, 8.$$ 
Observe also that $n_{1} = 1$, $n_{2} = 21$, $n_{3} = 56$, $n_{4} = 42$, 
$n_{7} = 48$ and $n_{j} = 0$ for other $j$ if $G = L_{2}(7)$, where 
$n_{d}$ denotes the cardinality of the elements of order $d$ in $G$.    
Now combining all of these together, we obtain
$$2+ \text{\rm rank}\, L^G = \frac{1}{168}(24 + 8 \times 21 + 6 \times 56 + 
4 \times 42 + 3 \times 48) = 5.
$$  
The remaining assertions now follow from the facts summarized before Claim 
(2.1).
\qed \enddemo
\remark{Remark (2.2)} By (2.1)(1), K3 surfaces with $L_{2}(7)$-action are 
of the maximum Picard number 20. By a similar case-by-case calculation, 
one can also show that the invariant lattices are positive definite and of 
rank $3$ if a K3 surface admits one of the eleven maximum symplectic group 
actions listed in [Mu1]. In particular, one has that 
\roster 
\item the invariant lattices (tensorized by $\bold R$) have the hyperk\"ahler 
three-space structure;
\item  such {\it algebraic} K3 surfaces are of the maximum Picard number 20 and are then at most countably many by [SI]. 
\endroster 
It would be very interesting to describe all of such (algebraic) K3 surfaces as rational points of the twister spaces corresponding to the invariant lattices (tensorized by $\bold R$). 
\qed \endremark 
Next we determine the discriminant of $L^{G}$.
\proclaim{ Key Lemma} $\vert \text{\rm det} L^{G} \vert = 196$.  \endproclaim
The proof of Key Lemma will be given after Claim (2.6).
Technically, this is the most crucial step and the next embedding Theorem 
due to Kondo is the most important ingredient in our proof of Key Lemma: 

\proclaim{Theorem (2.3) [Ko1]} Under the notation explained in Section 1, 
one has the following:
\roster
\item For a given finite symplectic action $H$ on a $K3$ surface,
there exists a non-Leech Niemeier lattice $N$ such that $L_{H} \subset N$.
Moreover, the action of $H$ extends to an action on $N$ so that
$L_{H} \simeq N_{H}$ and that $N^{H}$ contains a simple root.
\item This group action of $H$ on $N$ preserves the distinguished Weyl chamber 
$\Cal C$ and the natural homomorphism $H \rightarrow S(N)$ is injective. \qed
\endroster
\endproclaim
\proclaim{Corollary (2.4)} Under the notation of Theorem (2.3), one has: 
\roster 
\item $\text{\rm rank}\, N^{H} = \text{\rm rank}\, L^{H} + 2$. 
In particular, $\text{\rm rank}\, N^{G} = \text{\rm rank}\, L^{G} + 2 = 5$. 
\item $\vert \text{\rm det}\, N^{H} \vert = 
\vert \text{\rm det}\, L^{H} \vert$.  \qed
\endroster 
\endproclaim

\demo{Proof} Since $\text{rank}\, N^{H} = 24 - \text{rank}\, N_{H}$ 
and $\text{rank}\, L^{H} = 22 - \text{rank}\, L_{H}$, the first part of the 
assertion (1) follows from $N_{H} \simeq L_{H}$. Now the last part of (1) 
follows from (2.1). Recall that $L$ and $N$ are unimodular and 
the embeddings $L^{H} \subset L$ and $N^{H} \subset N$ are primitive. 
Then $\vert \text{det}L^{H} \vert = \vert \text{det}L_{H} \vert$ 
and $\vert \text{det}N^{H} \vert = \vert \text{det}N_{H} \vert$. 
Combining these with $N_{H} \simeq L_{H}$, one 
obtains $\vert \text{det}N^{H} \vert = 
\vert \text{det}L^{H} \vert$. \qed \enddemo  

Let us return back to our original situation and determine the Niemeier
lattice $N$ for our $G$. Note that $N$ is not the Leech lattice by (2.3)(1). 

\proclaim{Claim (2.5)} The Niemeier lattice $N$ in Theorem (2.3) for
$G = L_{2}(7)$ is $N(A_{1}^{\oplus 24})$.
\endproclaim

\demo{Proof} By Theorem (2.3)(2), $168 = |G|$ divides $|S(N)|$.
Thus, $N$ is either $N(A_{1}^{\oplus 24})$ or $N(A_{3}^{\oplus 8})$ 
by (1.3).
Suppose that the latter case occurs. Then, by (1.3), we have
$G \subset C_{2} \rtimes (C_{2}^{\oplus 3} \rtimes L_{3}(2))$.
Since $G$ is simple, a normal subgroup $G \cap  C_{2}$ is trivial, i.e.
$G \subset  C_{2}^{\oplus 3} \rtimes L_{3}(2)$. Again, for the same reason,
one has $G \subset L_{3}(2)$ (and in fact equal).
The Dynkin diagram $A_{3}^{\oplus 7}$, the complement of
the identity component in $A_{3}^{\oplus 8}$, consists of
the 21 simple roots $r_{i1}, r_{i2}, r_{i3}$ ($1 \leq i \leq 7$)
such that $(r_{i1}, r_{i2}) =
(r_{i2}, r_{i3}) = 1$ but $(r_{i1},r_{i3}) = 0$, $(r_{ik}, r_{jl}) = 0$ if $i
\not= j$.
Therefore the action $G$ on these  21 simple roots satisfies $g(r_{i2}) =
r_{g(i)2}$,
where $g$ in the right hand side is regarded as an element of
the permutation of the seven components.
Therefore $g(r_{i1})$ is either $r_{g(i)1}$ or $r_{g(i)3}$.
Thus, $G$ is embedded into a subgroup of the permutation subgroup
$C_{2}^{\oplus 7} \rtimes S_{7}$
of the 14 simple roots $r_{i1}, r_{i3}$. Here, the indices $1, 3$ are
so labelled that $\sigma \in S_{7}$ acts as $\sigma(r_{i1}) =
r_{\sigma (i) 1}$ and
$\sigma(r_{i3}) = r_{\sigma (i) 3}$  and the $i$-th factor of
$C_{2}^{\oplus 7}$
acts as a permutation of $r_{i1}$ and $r_{i3}$. Since $G$ is simple and
can not be embedded in $C_{2}^{\oplus 7}$, we have $G \subset S_{7}$.
Therefore,
$g(r_{i1}) = r_{g(i)1}$ and  $g(r_{i3}) = r_{g(i)3}$. In conclusion,
the orbits of the action $G$ on the 24 simple roots are $\{r_{01}\}$,
$\{r_{02}\}$,  $\{r_{03}\}$,  $\{r_{i1} \vert 1 \leq i \leq 7\}$,
$\{r_{i2} \vert 1 \leq i \leq 7\}$, $\{r_{i3} \vert 1 \leq i \leq 7\}$.
In particular, the 24 simple roots are divided into exactly
6 $G$-orbits.
Since these 24 roots generate the Niemeier lattice
$N = N(A_{3}^{\oplus 8})$ over $\bold Q$, we have then $\text{rank}\, N^{G} 
= 6$, a contradiction to (2.4)(1). Hence
the Niemeier lattice for our $G$ is $N(A_{1}^{\oplus 24})$.  \qed
\enddemo

From now we set $N := N(A_{1}^{\oplus 24})$. By (2.5) and (1.4), we have
$$G \subset M_{24} \subset S_{24} = \text{Aut}_{\text{set}}
(R).$$ 
Here $R := \{r_{i} \}_{i = 1}^{24}$ is the set of the simple roots 
of $N$ and the last inclusion is the natural one explained in (1.4). 
This allows us to use the table of the cyclic types of elements of 
$M_{24}$ given in [EDM] for its action on $R$. One may also talk about the 
orbit decomposition type of the action of $G$ on $R$. Although we donot know 
much about how $G$ is embedded 
in $M_{24}$, we can say at least the following: 

\proclaim{Claim (2.6) (cf. [Mu2])} The orbit decomposition type of $R$ by $G$ is either
$$[14, 1, 1, 7, 1] \,\, \text{or} \,\, [8, 7, 1, 7, 1].$$
\endproclaim

\demo{Proof} Since $N^{G} =  5$ by (2.4)(1), the 24 simple roots of $N$ are 
divided into exactly 5 $G$-orbits. (See the last argument of the Claim (2.5).) 
Set the orbit decomposition type as $[a, b, c, d, e]$.
Then $a + b + c+ d +e = 24$ and each entry is less than $21$. In addition,
since $G$
is simple and contains an element of order $7$, if $a \leq 6$ then $a = 1$,
for otherwise the natural non-trivial representation $G \rightarrow S_{a}$
would have a non-trivial kernel. Moreover, if $a \geq 7$, then $a$
divides $168 = \vert G \vert$.
This is because the action of $G$ on each orbit is, by the definition,
transitive.
Therefore $a$ is either $1$, $7$, $8$, $12$ or $14$. If $a = 12$,
then an element of order $7$ in $G$ has already $5$ fixed points
in this orbit. However, by [EDM], the cycle type of order $7$
element in $M_{24}$ is $(7)^{3}(1)^{3}$ and therefore has only $3$ fixed
points,
a contradiction. Hence $a$ is either $1$, $7$, $8$ or $14$. Clearly
the same holds for $b$, $c$, $d$, $e$. Now by combining these, together with
the equality $a + b + c+ d +e = 24$, we obtain the result. \qed
\enddemo

\demo{Proof of Key Lemma} By (2.4)(2), we may calculate 
$\vert \text{det}N^{G} \vert $ instead.
Let us renumber the 24 simple roots according to the orbit decompositions
found in (2.6): 
$$\{r_{1}, \cdots , r_{14}\} \cup \{r_{15}\} \cup \{r_{16}\}
\cup \{r_{17}, \cdots ,  r_{23}\} \cup \{r_{24}\} - (*)$$ 
or 
$$\{r_{1}, \cdots , r_{8}\} \cup \{r_{9}, \cdots , r_{15} \}
\cup \{r_{16}\} \cup \{r_{17}, \cdots ,  r_{23}\} \cup \{r_{24}\} -(**).$$
\par \vskip 4pt
Consider the case ($*$) first.
Recall that rank $N^{G} = 5$ and $N^{G} = N \cap (A_{1}^{\oplus 24})^{G} 
\otimes \bold Q$. 
Then
$b_{1} = \sum_{i=1}^{14}r_{i}$, $b_{2} = r_{15}$,
$b_{3} = r_{16}$, $b_{4} = \sum_{i=17}^{23}r_{i}$ and  $b_{5} =  r_{24}$
form the basis of $(A_{1}^{\oplus 24})^{G}$. Moreover, by (1.4), we see 
that $N^{G}/(A_{1}^{\oplus 24})^{G}$ consists of the elements 
of the form $\sum_{i \in I} b_{i} / 2$, where the set of the simple roots 
$\{r_{j}\}$ appearing in the sum $\sum_{i \in I} b_{i}/2$
is either $R$, $\emptyset$, an Octad, complement of an Octad, or a Dodecad.
However, by the shape of the orbit decomposition, there are no cases where 
a Dodecad appears.
Therefore, in order to get an integral basis of $N^{G}$ we may find out
all the Octads and their complements appearing in the forms above.

\proclaim{Claim (2.7)} By reordering the three 1-element orbits 
if necessary, the union of the fourth and fifth orbits $ \{r_{17}, 
r_{18}, \cdots , r_{23}, r_{24} \}$ forms an Octad. \endproclaim
\demo{Proof}  Let $\alpha \in G$ be an element 
of order 7. Then the cycle type of $\alpha$ (on $R$) is $(7)^{3}(1)^{3}$ [EDM]. 
In particular, a simple root $x$ 
forms a 1-element orbit if $\alpha^{k}(x) = x$ for some $k$ 
with $1 \leq k \leq 6$. Moreover, one can adjust the numbering of the roots 
in the fourth orbit $\{r_{17}, r_{18}, \cdots , r_{23}\}$ so as to 
be that $\alpha(r_{i}) = r_{i+1}$ ($17 \leq i \leq 22$) and $\alpha(r_{23}) = 
r_{17}$. Let us consider the 5-element set $S := 
\{r_{17}, r_{18}, \cdots , r_{21} \}$. Then by (1.5), there is an Octad $O$ 
such that 
$S \subset O$. We shall show that $O = \{r_{17}, 
r_{18}, \cdots , r_{23}, r_{24} \}$. For this purpose, assuming first that 
$r_{22}, r_{23} \not\in O$, we shall derive a contradiction. Under this 
assumption, one has 
$\{r_{19}, r_{20}, r_{21}\} \subset O \cap \alpha^{2}(O)$ and $(O \cap 
\alpha^{2}(O)) \cap \{r_{17}, r_{18}, r_{22}, r_{23} \} 
= \emptyset$. (Here for the last equality we used the fact that 
$\alpha^{2}(r_{22}) = r_{17}$ and 
$\alpha^{2}(r_{23}) = r_{18}$.)
The last equality also implies that $O \not= \alpha^{2}(O)$.  
Let us consider the symmetric difference 
$D = (O - \alpha^{2}(O)) \cup (\alpha^{2}(O) - O)$. 
Then $\vert D \vert$ must be either $8$, $12$, $16$ whence 
$\vert O \cap \alpha^{2}(O) \vert = 4$ (See the proof 
of (1.5)).  
Thus, 
there exists an $x \in R$ such that 
$O \cap \alpha^{2}(O) = \{r_{19}, r_{20}, r_{21}, x \}$ and one has 
$$O = \{r_{17}, r_{18}, r_{19}, r_{20}, r_{21}, x, y, z \}.$$ 
Here none of $x, y, z$ lies in the fourth orbit.  
Since $x \in \alpha^{2}(O)$, one has either $x = \alpha^{2}(y)$ or 
$x = \alpha^{2}(x)$ (by changing the role of $y$ and $z$ if neccesary). 
In each case, we have $\alpha^{2}(z) \not= z$, whence $\alpha^{k}(z) \not= z$ 
for all $k$ with $1 \leq k \leq 6$. In particular, $z$ is in the first orbit. 
\par 
Consider first the case where $x = \alpha^{2}(y)$. In this case, both $x$ 
and $y$ 
belong to the first orbit. 
Let us rename the elements in the first orbit so as to be that 
$\alpha(r_{i}) = r_{i+1}$ 
($1 \leq i \leq 6$), $\alpha(r_{7}) = r_{1}$; $\alpha(r_{7+i}) = r_{i+8}$ 
($1 \leq i \leq 6$), $\alpha(r_{14}) = r_{8}$ and that $y = r_{1}$ and 
$x = r_{3}$. Then, we have 
$$O = \{r_{17}, r_{18}, r_{19}, r_{20}, r_{21}, r_{3}, r_{1}, z \},$$
and 
$$\alpha^{3}(O) = \{r_{20}, r_{21}, r_{22}, r_{23}, r_{17}, r_{6}, r_{4}, 
\alpha^{3}(z) \}.$$ 
Considering the symmetric difference of $O$ and $\alpha^{3}(O)$ as before, 
one finds that $O \cap \alpha^{3}(O)$ is a 4-element set. 
Therefore $\vert \{r_{1}, r_{3}, z\} \cap \{r_{4}, r_{6}, \alpha^{3}(z) \} 
\vert = 1$. Combining this with $\alpha^{3}(z) \not= z$, one has either 
$r_{1} = \alpha^{3}(z)$, $r_{3} = \alpha^{3}(z)$, $z = r_{4}$, or $z = r_{6}$. 
Hence $O$ satisfies 
one of the following four:
$$O = \{r_{17}, r_{18}, r_{19}, r_{20}, r_{21}, r_{1}, r_{3}, r_{5} \} - (1)$$
$$O = \{r_{17}, r_{18}, r_{19}, r_{20}, r_{21}, r_{1}, r_{3}, r_{7} \} - (2)$$
$$O = \{r_{17}, r_{18}, r_{19}, r_{20}, r_{21}, r_{1}, r_{3}, r_{4} \} - (3)$$
$$O = \{r_{17}, r_{18}, r_{19}, r_{20}, r_{21}, r_{1}, r_{3}, r_{6} \} - (4).$$ In the case (1), one calculates 
$$\alpha^{2}(O) =  \{r_{19}, r_{20}, r_{21}, r_{22}, r_{23}, r_{3}, r_{5}, r_{7} \}$$
and has then $O \cap \alpha^{2}(O) =  \{r_{19}, r_{20}, r_{21}, r_{3}, 
r_{5}\}$. In particular, the two Octads $O$ and $\alpha^{2}(O)$ share 5 
elements in common. Then, by the Steiner property (1.5), we would have 
$O = \alpha^{2}(O)$, a contradiction. By considering 
$O \cap \alpha(O)$ in the cases (2), (3) and $O \cap \alpha^{2}(O)$ in the case (4), we can derive a contradiction in the same manner, too. Thus, the case 
$x = \alpha^{2}(y)$ is impposible. 
\par 
Next we consider the case where $\alpha^{2}(x) = x$. In this case,  
this $x$ forms a 1-element orbit and satisfies $\{r_{18}, r_{19}, r_{20}, 
r_{21}, x \} \subset O \cap \alpha(O)$. 
However, the Steiner property would then 
imply $O = \alpha(O)$, whence $O = \alpha^{2}(O)$, a contradiction. 
\par
Therefore, 
the Octad $O$ satisfies either $r_{22} \in O$ or $r_{23} \in O$, i.e.  
$$\{r_{17}, r_{18}, r_{19}, r_{20}, r_{21}, r_{22}\} \subset O,$$ 
or 
$$\{r_{23}, r_{17}, r_{18}, r_{19}, r_{20}, r_{21}\} \subset O.$$  
Then one has either 
$$\{r_{18}, r_{19}, r_{20}, r_{21}, r_{22}\} \subset O \cap \alpha(O),$$ 
or 
$$\{r_{17}, r_{18}, r_{19}, r_{20}, r_{21}\} \subset O \cap \alpha(O).$$
Hence, by the Steiner property, we have $O = \alpha(O)$, whence 
$O = \alpha^{k}(O)$ for all $k$. This implies that the Octad $O$ is of the 
form  
$$O = \{r_{17}, r_{18}, r_{19}, r_{20}, r_{21}, r_{22}, r_{23}, x\}$$ 
for some root $x$. Since $O = \alpha(O)$, we have also $\alpha(x) = x$. 
Hence this $x$ forms a 1-element orbit set. 
\qed \enddemo              
By this Claim, one has $b_{6} :=   (b_{4} + b_{5})/ 2 \in N^{G}$ and also
$b_{7} := (b_{1} + b_{2} + b_{3})/ 2 \in N^{G}$.
By the remark before Claim (2.7) and the Steiner property (1.5), we also see 
that 
there are no other Octads
appearing in the sum $\sum_{i \in I} b_{i}/2$.
Since $\sum_{i=1}^{5}b_{i}/2 = b_{6} + b_{7}$, the seven elements
$b_{1}, \cdots, b_{7}$ then generate $N^{G}$ over $\bold Z$.
Moreover, since $b_{1} = 2b_{7} - b_{2} - b_{3}$ and $b_{4} = 2b_{6} - b_{5}$,
we finally
see that $b_{7}, b_{2}, b_{3}, b_{6}, b_{5}$ form an integral basis of
$N^{G}$.  Using $(r_{i}, r_{j}) = -2\delta_{ij}$, we
find that the intersection matrix
of $N^{G}$ under this basis is given as $A$ below:
$$A = \pmatrix -8 & -1 & -1 &  0 &  0 \\
               -1 & -2 &  0 &  0 &  0 \\
               -1 &  0 & -2 &  0 &  0 \\
                0 &  0 &  0 & -4 & -1 \\
                0 &  0 &  0 & -1 & -2
\endpmatrix, \,\,\,\,
B = \pmatrix -4 &  0 &  0 &  0 &  0 \\
              0 & -4 & -1 &  0 &  0 \\
              0 & -1 & -2 &  0 &  0 \\
              0 &  0 &  0 & -4 & -1 \\
              0 &  0 &  0 & -1 & -2
\endpmatrix.$$
\par 
Let us consider next the case ($**$). Since $G$ acts on the first 8-element 
orbit 
transitively, for each root $r_{i}$ in the first orbit, one can find an element $\alpha' \in G$ of order $7$ such that $\alpha'(r_{i}) \not= r_{i}$. Now, 
by the same argument 
based on the fact that the cycle type of order 7 element is 
$(7)^{3}(1)^{3}$ and the Steiner property (1.5) 
(together with the remark above), 
one finds that (after reordering the two 1-element orbits) the union 
of the second and third orbits and the union of the fourth and fifth 
orbits are both Octads. This also implies that the first orbit is an Octad. 
Then again by the same argument as in the previous case, one can easily 
see that the five elements
$b_{1} := \sum_{i = 1}^{8} r_{i}/2$, $b_{2} := \sum_{i = 9}^{16} r_{i}/2$,
$b_{3} = r_{16}$, $b_{4} = \sum_{i = 17}^{24} r_{i}/2$, and $b_{5} = r_{24}$
form an integral basis of $N^{G}$ in the second case, and that 
the intersection matrix of $N^{G}$ under this basis is given as $B$ above.
Clearly $|\text{det} N^{G}| = 196$ in both cases.
This proves Key Lemma. \qed \enddemo
Next we shall study possible extensions
$G \subset \tilde{G}$. Recall that we have already shown that $\tilde{G}/G 
\simeq \mu_{I}$,
where $I$ is either $1$, $2$, $3$, $4$ or $6$, and that
$\tilde{G}/G$ acts faithfully on
$T_{X}$. Set $\tilde{G}/G = \langle \tau \rangle$.
\par
\vskip 4pt
The next Lemma is valid for any $G \subset \tilde{G}$ if
$\tilde{G}/G$ acts faithfully on $T_{X}$ and if rank $T_{X} = 2$.

\proclaim{Lemma (2.8)}
\roster
\item Assume that $\text{ord}(\tau) = 3$. Then,
as $\bold Z[\tau^{*}]$-modules,
one has $T_{X} \simeq \bold Z[x]/(x^{2}+x+1)$, where $\tau^{*}$ acts on the
right
hand side by the multiplication by $x$. In particular, one can take an integral basis
$e_{1}, e_{2}$ of  $T_{X}$ such that
$\tau^{*}(e_{1}) =  e_{2}$ and $\tau^{*}(e_{2}) =  -(e_{1} + e_{2})$.
Moreover, under this basis, the intersection matrix of $T$ is of
the form $((e_{i},e_{j})) =$
$\pmatrix 2m & -m \\
          -m & 2m
\endpmatrix.$
\item Assume that $\text{ord}(\tau) = 4$. Then, as $\bold Z[\tau^{*}]$-modules,
one has
$T _{X}\simeq \bold Z[x]/(x^{2}+1)$, where $\tau^{*}$ again
acts on the right hand side by
the multiplication by $x$. In particular, one can take an integral basis
$e_{1}, e_{2}$ of
$T_{X}$ such that $\tau^{*}(e_{1}) = e_{2}$ and $\tau^{*}(e_{2}) = -e_{1}$.
Moreover, under this basis, the intersection matrix of $T$ is of the form
$((e_{i},e_{j})) =$
$\pmatrix 2m &  0 \\
           0 & 2m
\endpmatrix.$
\endroster
\endproclaim
\demo{Proof} The first part of the two assertions is due to the fact that
$\bold Z[\zeta_{3}]$ and $\bold Z[\zeta_{4}]$ are both PID.
(For more detail, see for example [MO].) By taking
an integral basis
of $T_{X}$ corresponding to $1$ and $x$ (in the right hand side),
one obtains the desired representation of the action of $\tau^{*}$.
Now combining this with $(\tau^{*}(a),\tau^{*}(b)) = (a,b)$,
we get the intersection matrix as claimed. \qed \enddemo

The next Claim completes the first assertion of the main Theorem:

\proclaim{Claim (2.9)} $I \not= 6$. \endproclaim 

A similar method is exploited in [Ko2] and [OZ] in other settings with 
somewhat different flavours and will be also adopted in the next Claim (2.10).

\demo{Proof} Assuming to the contrary that $\tilde{G}/G = \langle g \rangle 
\simeq \mu_{6}$, we shall 
derive a contradiction. By (2.1), one has $L^{G} \supset T_{X} 
\oplus \bold Z H$, where
$H$ is the primitive ample class invariant under $G$. Set $(H^{2}) = 2n$.
Since $T_{X}$ is primitive in $L^{G}$,
we can choose an integral basis of
$L^{G}$ as $e_{1}, e_{2}$ and $e_{3} = (aH + be_{1} + ce_{2})/\ell$,
where $e_{1}$ and $e_{2}$ are  the integral basis of $T_{X}$ found
in (2.8)(1) applied for $\tau := g^{2}$ and $\ell$ and $a, b, c$ are 
integers 
such that $(\ell, a) = 1$.
Then 
$$L^{G}/ (T_{X} \oplus \bold Z H)  =
\langle \overline{e_{3}} \rangle \simeq C_{\ell},$$
where $ \overline{e_{3}} = e_{3}$ mod ($T_{X} \oplus \bold Z H$). Since $H$
is also stable under $\tilde{G}$, we have $\tau^{*}(\overline{e_{3}}) =
\overline{e_{3}}$
and 
$$\tau^{*}(be_{1} + ce_{2})/\ell \equiv (be_{1} + ce_{2})/\ell \,  \text{mod} 
\, T_{X}.$$
On the other hand, by the choice of $e_{1}, e_{2}$, we calculate
$$\tau^{*}(be_{1} + ce_{2})/\ell = (-ce_{1} + (b - c)e_{2})/\ell.$$ 
Therefore,
$b \equiv -c$ and
$c \equiv b - c$ $\text{mod}\, \ell$. In particular,
$b \equiv -c$ and $3b \equiv 3c \equiv 0$
$\text{mod}\, \ell$. This, together with the primitivity of
$\bold Z H$ in $L^{G}$, implies that $\ell = 1$
or $3$, that is, $[L^{G} :  T_{X} \oplus \bold Z H] = 1$ or $3$. 
\par
\vskip 4pt
If $\ell = 1$, we have $L^{G} =  T_{X} \oplus \bold Z H $ and $196 = 6m^{2}n$.
However $6$ is not a divisor of $196$, a contradiction. 
\par
\vskip 4pt 
Consider the case $\ell = 3$. Then, (by using the primitivity 
of $H$ in $L^{G}$ and by adding an element of $T_{X}$ to $e_{3}$ if 
necessary), we can take one of $(\pm H \pm (e_{2} - e_{3}))/3$ as $e_{3}$. 
Put $\sigma := g^{3}$. Then $\sigma^{*}H = H$ and $\sigma^{*} \vert T_{X} 
= -id$. Using these two equalities, we calculate 
$$\sigma^{*}(e_{3}) = \sigma^{*}((\pm H \pm (e_{2} - e_{3}))/3) = 
(\pm H \mp (e_{2} - e_{3}))/3.
$$ 
However, one would then have  
$$\pm 2(e_{2} - e_{3})/3 = e_{3} - \sigma^{*}(e_{3}) \in L^{G},$$ 
a contradiction to the primitivity of $T_{X}$ in $L^{G}$. 
Hence $I \not= 6$.   
\qed \enddemo
From now, we consider the maximum case $\tilde{G}/G = \langle \tau \rangle 
\simeq \mu_{4}$.  
\proclaim{Claim (2.10)} $(H^{2}) = 4$. \endproclaim 

\demo{Proof} As in (2.9), one has $L^{G} \supset T_{X} \oplus \bold Z H$, 
where$H$ is the primitive ample class invariant under $G$. Set $(H^{2}) = 2n$.
Since $T_{X}$ is primitive in $L^{G}$,
we can choose an integral basis of
$L^{G}$ as $e_{1}, e_{2}$ and $e_{3} = (aH + be_{1} + ce_{2})/\ell$,
where $e_{1}$ and $e_{2}$ are  the integral basis of $T_{X}$ found
in (2.8)(2) and $\ell$ and $a, b, c$ are integers such that $(\ell, a)
= 1$.
Then, as in (2.9), we have  
$$L^{G}/ (T_{X} \oplus \bold Z H)  =
\langle \overline{e_{3}} \rangle \simeq C_{\ell},$$
where $ \overline{e_{3}} = e_{3}$ mod ($T_{X} \oplus \bold Z H$). Since $H$
is also stable under $\tilde{G}$, we have $\tau^{*}(\overline{e_{3}}) =
\overline{e_{3}}$
and 
$$\tau^{*}(be_{1} + ce_{2})/\ell \equiv (be_{1} + ce_{2})/\ell \, 
\text{mod} \, T_{X}.$$
On the other hand, by the choice of $e_{1}, e_{2}$, we calculate
$$\tau^{*}(be_{1} + ce_{2})/\ell = (be_{2} -  ce_{1})/\ell.$$ 
Therefore,
$b \equiv c$ and
$c \equiv -b$ $\text{mod}\, \ell$. In particular,
$b \equiv c$ and $2b \equiv 2c \equiv 0$
$\text{mod}\, \ell$. This, together with the primitivity of
$\bold Z H$ in $L^{G}$, implies that $\ell = 1$
or $2$, that is, $[L^{G} :  T_{X} \oplus \bold Z H ] = 1$ or $2$. 
\par 
In the first case,
we have $ L^{G} =  T_{X} \oplus \bold Z H $ and $196 = 8m^{2}n$.
However $8$ is not a divisor of $196$, a contradiction. 
\par 
In the second case,
we have $2^{2} \cdot 196 = 8m^{2}n$, i.e. $m^{2}n = 2 \cdot 7^{2}$. Then
$(m,n)$ is either
$(1, 2 \cdot 7^{2})$ or $(7, 2)$. In the first case we have $X = X_{4}$ by
the result of Shioda and Inose [SI],
where $X_4$ is the minimal resolution of
$(E_{\sqrt{-1}} \times E_{\sqrt{-1}})/ \langle
\text{\rm diag}(\sqrt{-1}, -\sqrt{-1}) \rangle $.
However, according to the explicit description of $\text{Aut}(X_{4})$ by
Vinberg [Vi],
$X_{4}$ has no automorphism of order $7$, a contradiction.
Therefore, only the second case can happen and one has $(H^{2}) = 2n = 4$ 
(and $T_{X} = \text{diag}(14, 14)$).
\qed \enddemo

Now the following Claim will complete the proof of the main Theorem.

\proclaim{Claim (2.11)} $(X, \tilde{G})$ is isomorphic to
$(X_{168}, L_{2}(7) \times \mu_{4})$
defined in the Introduction.
\endproclaim

\demo{Proof} Since $S_{X}^{G} = \bold Z H$,
$\vert H \vert$ has no fixed components.
Indeed, the fixed part of $\vert H \vert$ must be also $G$-stable but
is of negative definite [SD].
Therefore, the ample linear system $\vert H \vert$ is free 
[ibid.]. Note that  $\text{dim} \vert H \vert = 3$ by
the Riemann-Roch formula and the fact $(H^{2}) = 4$. Then $\vert H \vert$
defines a morphism $\Phi := \Phi_{\vert H \vert} : X \rightarrow \bold P^{3}$.
This $\Phi$ is either an embedding to a quartic surface $S$ or a finite double
cover of
an integral quadratic surface $Q$.  Note that $\tilde{G}$ acts on the image
as a projectively linear transformation. Moreover, the action of $G$ on the image is 
faithful
even in the second case, because $G$ is simple.
Recall that the degrees of the projectively linear irreducible representations 
of
$G$ are $1, 3, 4, 6, 7, 8$ [ATLAS] and that the two 3-dimensional irreducible
representations are transformed by the outer automorphism of $G$.
Then the action of $G$ on the image is induced by the irreducible decomposition
${\bold C}^4 = V_{1} \oplus V_{3}$ or ${\bold C}^{4} = V_{4}$. (Note that the action of $G$ on $\bold P^{3}$ is linearlized in the first case but not in the second case. More precisely, in the second case, only the action of $\text{SL}(2, \bold F_{7})$, i.e. the central extension of $G$ corresponding to the Schur multiplier $2$, is linearlized.) 
\par
\vskip 4pt 
Let us first consider the second case. By [Ed, Pages 198 - 200 and Page 166], $G$ has no invariant hypersurface of degree $2$ but only one invariant hypersurface of degree $4$: 
$$f := x_{0}^{4} + 6\sqrt{2}x_{0}x_{1}x_{2}x_{3} + x_{1}x_{2}^3
+ x_{2}x_{3}^3 + x_{3}x_{1}^3 = 0,$$
where the homogeneous coordinates $[x_{0} : x_{1} : x_{2} : x_{3}]$ are 
chosen in such a way that an order $7$-element of $G$, say $\alpha$, is represented by the 
following diagonal matrix: 
$$A = \pmatrix  1 &  0 &  0 &  0 \\
                0 &  \zeta_{7}^{6} &  0 &  0 \\
                0 &  0 &  \zeta_{7}^{3} &  0 \\
                0 &  0 &  0 & \zeta_{7}^{5}\endpmatrix .$$ 
Thus, $\Phi$ is an embedding and one has $S = (f = 0)$. Recall that 
$\tilde{G}$ fits in with the exact sequence 
$$ 1 \rightarrow G \rightarrow \tilde{G} \rightarrow \mu_{4} \rightarrow 1,$$
where the last map is the representation of $\tilde{G}$ on $H^{0}(S, \Omega_{S}^{2}) = \bold C \omega_{S}$. Let $g \in \tilde{G}$ be a lift of $\zeta_{4} \in \mu_{4}$. Then $g^{*}\omega_{S} = \zeta_{4}\omega_{S}$. Since $G$ is a normal subgroup of $\tilde{G}$, one can define an element $c_{g} \in \text{Aut}_{\text{group}}(G)$ by $G \ni x \mapsto g^{-1}xg \in G$ for all $x \in G$. Since $\text{Out}(G) \simeq C_{2}$ [ATLAS], $(c_{g})^{2}$ is then an inner automorphism of $G$, i.e. there exists an element $y \in G$ such that $g^{-2}xg^{2} = y^{-1}xy$ for all $x \in G$. Set $k = g^{2}y^{-1}$. Then one has 
$k^{-1}xk = x$ for all $x \in G$, $k^{*}\omega_{S} = -\omega_{S}$ and $2 \vert \text{ord}(k)$. 
Therefore, replacing $k$ by $k^{2l+1}$ if necessary, one obtains an element 
$h \in \tilde{G}$ such that $h^{-1}xh = x$ for all $x \in G$, $h^{*}\omega_{S} = -\omega_{S}$ and $\text{ord}(h) = 2^{n}$. Choose a representative $(h_{ij})$ of $h$ in $\text{GL}(4, \bold C)$. Then for $A$ above one has 
$(h_{ij})A(h_{ij})^{-1} = cA$ ($c \in \bold C$) in $\text{GL}(4, \bold C)$. 
This implies $c^{2^{n}} = 1$ and $\{1, \zeta_{7}^{6}, \zeta_{7}^{3}, 
\zeta_{7}^{5} \} = \{c, c\zeta_{7}^{6}, c\zeta_{7}^{3}, 
c\zeta_{7}^{5} \}$. Thus $c = 1$ and one has $(h_{ij})A(h_{ij})^{-1} = A$, 
i.e. $(h_{ij})A = A(h_{ij})$ in $\text{GL}(4, \bold C)$. This readily implies that $(h_{ij})$ is also a diagonal matrix, and one may write that 
$$(h_{ij}) = \pmatrix  1 &  0 &  0 &  0 \\
                0 &  h_{11} &  0 &  0 \\
                0 &  0 &  h_{22} &  0 \\
                0 &  0 &  0 & h_{33}\endpmatrix .$$ 
Then, by the shape of $f$ and by the fact that $h^{*}f = c'f$ for some $c' \in \bold C$, one has $c' = 1$ and $h_{11}h_{22}h_{33} = 1$. However, this would yield $(h_{ij})^{*}f = f$ and $\text{det}(h_{ij}) = 1$, thereby $h^{*}\omega_{S} = \omega_{S}$, a contradiction to the previous equality $h^{*}\omega_{S} = -\omega_{S}$. Hence the second case cannot happen. 
\par 
\vskip 4pt 
Let us next consider the first case.   
Let us choose the homogeneous
coordinates $[x_{0} : x_{1} :  x_{2} : x_{3} ]$ such that $x_{0}$ is
the coordinate of $V_{1}$ and that $x_{i}$ ($1 \leq i \leq 3$) are the 
coordinates of $V_{3}$ described as in the Introduction. 
\par
\vskip 4pt
Let us first consider the case where $\Phi$ is a double covering.
Write an equation of $Q$ as
$$ax_{0}^{2} + x_{0}f_{1}(x_{1}, x_{2}, x_{3}) +  f_{2}(x_{1}, x_{2}, x_{3})
= 0.$$
Since $G$ is simple and acts on $x_{0}$ as an identity, we have $g^{*}(f_{1})
= f_{1}$
and $g^{*}(f_{2}) = f_{2}$ for all $g \in G$.
Since there are no non-trivial $G$-invariant
linear and quadratic forms in three variables [Bu, Section 267], 
one has $f_{1} = f_{2} = 0$ .
However, then $Q $ is not integral, a contradiction. Hence $\Phi$ is
an embedding. Let us write an equation of $S$ as 
$$F =
ax_{0}^{4} + x_{0}^{3}f_{1}(x_{1}, x_{2}, x_{3}) +
x_{0}^{2}f_{2}(x_{1}, x_{2}, x_{3}) +
x_{0}f_{3}(x_{1}, x_{2}, x_{3}) + f_{4}(x_{1}, x_{2}, x_{3}) = 0.$$
Then $g^{*}(f_{i}) = f_{i}$ for all $g \in G$ and all $1 \leq i \leq 4$.
Thus by [ibid.], we have $f_{i} = 0$ for all $1 \leq i \leq 3$,
$ f_{4}(x_{1}, x_{2}, x_{3}) = b( x_{1}x_{2}^3 + x_{2}x_{3}^3 + x_{3}x_{1}^3)$
and $F(x_{0}, x_{1}, x_{2}, x_{3}) =
ax_{0}^{4} +  b( x_{1}x_{2}^3 + x_{2}x_{3}^3 + x_{3}x_{1}^3)$.
Here $a \not= 0$ and $b \not= 0$, because $S$ is non-singular.
Therefore, by multiplying coordinates suitably, one may adjust the equation
of $S$ as 
$$x_{0}^{4} +  x_{1}x_{2}^3
+ x_{2}x_{3}^3 + x_{3}x_{1}^3 = 0.$$ 
Hence $S \simeq X_{168}$ and $L_{2}(7)
\times \mu_{4}$
acts on $S$ as described in the Introduction, where by the construction, the
action $L_{2}(7)$
also coincides with the given action of $G$ on $X$.
Since $S$ is a K3 surface and has no non-zero global holomorphic vector
fields,
the projectively linear automorphism group $G''$ of $S \subset {\bold P}^3$ is finite.
This $G''$ satisfies $G'' \supset L_{2}(7) \times \mu_{4}$ and 
$G'' \supset \tilde{G}$. Thus $4 \cdot 168 \vert \vert G'' \vert$ and 
one has $\vert G'' \vert = \vert  L_{2}(7) \times \mu_{4} \vert  =
\vert \tilde{G} \vert = 4 \cdot 168$.
Hence $\tilde{G} = G'' =  L_{2}(7) \times \mu_{4}$
as projectively linear automorphism groups of $S$. Now we are done. \qed \enddemo
\par \vskip 2pc 
\remark{Remark (2.12)}
\roster
\item By the proof of (2.10), we have that
$T_{X} =$ diag$(14, 14)$ if $\vert \tilde{G}/G \vert = 4$.
In particular, $T_{X_{168}} =$ diag$(14, 14)$.
\item Now one can easily check that the two K3 surfaces
$X_{168}$ and $X_{168}'$ in the Introduction are
not isomorphic to each other. Note that $X_{168}'$ has a $G$-stable ample 
class $H$ of degree $2$. Therefore, if $T_{X_{168}'}$ 
is isomorphic to $T_{X_{168}} =$ diag $(14, 14)$, 
then $[L^{G} : T_{X_{168}'} \oplus \bold{Z}H]^{2} = 
2 \cdot 14^{2}/196 = 2$, a contradiction.  \qed
\endroster
\endremark
\head References \endhead

\par \noindent
[ATLAS] J. H. Conway et al. ATLAS of Finite Groups, Clarendon Press - Oxford 
(1985).

\par \noindent
[Bu] W. Burnside, Theory of groups of finite order,
Dover Publications, Inc., New York (1955).

\par \noindent
[CS] J. H. Conway and N. J. A. Sloane, Sphere packings, lattices and groups,
Grundlehren der Math. Wissenschaften,
Springer-Verlag, New York (1999).

\par \noindent
[Ed] W. L. Edge, The Klein group in three dimensions,
Acta Math. 79 (1947) 153--222.

\par \noindent
[EDM] Encyclopedic Dictionary of Mathematics, Vol. II (Appendix B, Table 5.I),
Translated from the second Japanese edition,
MIT Press, Cambridge, Mass.-London (1977) pp. 885--1750.

\par \noindent
[Ko1] S. Kondo, Niemeier Lattices, Mathieu groups, and finite groups
of symplectic automorphisms of $K3$ surfaces,
Duke Math. J. {\bf 92} (1998) 593--598.

\par \noindent
[Ko2] S. Kondo, The maximum order of finite groups of automorphisms
of $K3$ surfaces, Amer. J. Math. {\bf 121} (1999) 1245--1252.

\par \noindent
[MO] N. Machida and K. Oguiso, On $K3$ surfaces admitting finite
non-symplectic group actions, J. Math. Sci. Univ. Tokyo {\bf 5} (1998) 
273--297.

\par \noindent
[Ma] D. Markushevich, Resolution of $\bold{C}^{3}/H_{168}$, 
Math. Ann. {\bf 308} (1997) 279--289.

\par \noindent
[Mu1] S. Mukai, Finite groups of automorphisms of $K3$ surfaces
and the Mathieu group, Invent. Math. {\bf 94} (1988) 183--221.

\par \noindent
[Mu2] S. Mukai, Lattice theoretic construction of symplectic
actions on $K3$ surfaces (an Appendix to [Ko1]),
Duke Math. J. {\bf 92} (1998)  599--603.

\par \noindent
[Ni] V. V. Nikulin, Finite automorphism groups of Kahler $K3$ surfaces,
Trans. Moscow Math. Soc. {\bf 38} (1980) 71--135.

\par \noindent
[Og] K. Oguiso, On the complete classification of Calabi-Yau threefolds
of Type $III_0$, in : Higher dimensional complex varieties,
de Gruyter (1996) 329--340.

\par \noindent
[OZ] K. Oguiso and D. -Q. Zhang, On the most algebraic $K3$ surfaces and
the most extremal log Enriques surfaces, Amer. J. Math. {\bf 118} (1996)
1277--1297.

\par \noindent
[SI] T. Shioda and H. Inose, On singular $K3$ surfaces,
in : Complex analysis and algebraic geometry, pp. 119--136,
Iwanami Shoten, Tokyo (1977).

\par \noindent
[SD] B. Saint-Donat, Projective models of $K3$ surfaces,
Amer. J. Math. {\bf 96} (1974) 602--639.

\par \noindent
[Vi] E. B. Vinberg, The two most algebraic $K3$ surfaces, 
Math. Ann. {\bf 265} (1983) 1--21.

\par \vskip 2pc
\head Added in Proof
\endhead

In this added in proof, we continue to employ the same notation, eg. $G = L_{2}(7)$. After submitting our note, we noticed the following Propositions. These answer questions asked by I. Dolgachev around January 2001: 
\proclaim{Proposition 1} In the main Theorem (1), one has $\vert \tilde{G}/G \vert \not= 3$. 
\endproclaim 
\proclaim{Proposition 2} There are infinitely many non-isomorphic algebraic K3 surfaces $X$ such that $G \subset \text{Aut}(X)$. 
\endproclaim 
In Proposition 2, we already know that there are at most countably many such algebraic K3 surfaces $X$ (Remark (2.2)). 
\head Proof of Proposition 1 \endhead 
Assuming to the contrary, we shall derive a contradiction. Recall that $\text{rk}S_{X} = 20$ and $\text{rk}T_{X} = 2$ if $X$ is an algebraic K3 surface admitting an $L_{2}(7)$-action.
\proclaim{Claim 1} $\tilde{G} \simeq G \times \mu_{3}$. \endproclaim 
\demo{Proof} Let $\tilde{g} \in \tilde{G}$ be an element such that 
$\tilde{g}^{*}\omega_{X} = \zeta_{3}\omega_{X}$. Replacing $\tilde{g}$
by its power $\tilde{g}^{n}$ such that $(n, 3) = 1$, one may assume that 
$\text{ord}(\tilde{g}) = 3^{m}$ for some positive integer $m$. Since $G$ is a 
normal subgroup of $\tilde{G}$, one has $c_{\tilde{g}}(x) := \tilde{g}^{-1}x \tilde{g} \in G$ if $x \in G$, thereby $c_{\tilde{g}} \in \text{Aut}(G)$. 
Since $\text{Out}(G) = C_{2}$ and $\text{ord}(\tilde{g}) = 3^{m}$, one has then $c_{\tilde{g}} \in \text{Inn}(G)$, i.e. there is $y \in G$ such that $\tilde{g}^{-1}x \tilde{g} = y^{-1}xy$ for all $x \in G$. Now, replacing $\tilde{g}$ by 
$\tilde{g}y^{-1}$, one has 
$\tilde{g}^{-1}x\tilde{g} = x$ for all $x \in G$ and 
$\tilde{g}^{*}\omega_{X} = \zeta_{3}\omega_{X}$. 
Then $\tilde{g}^{3} \in G$ and is also in the center of $G$. Note that 
the center of $G$ is $\{id.\}$ for $G$ being simple, non-commutative. Then 
$\tilde{g}^{3} = id.$ and $\tilde{g}$ gives a desired splitting of the exact sequence $1 \rightarrow G \rightarrow \tilde{G} \rightarrow \mu_{3} \rightarrow 1$. \qed 
\enddemo 
\proclaim{Claim 2} $T_{X}$ is isomorphic to $\pmatrix 14 & -7 \\
          -7 & 14
\endpmatrix$ (and the degree of the primitive invariant polarization is 
$12$). \endproclaim
\demo{Proof} The argument is the same as in (2.10) but is based on the following three facts instead: Lemma (2.8)(1) (instead of (2)); 
$[L^{G} :  T_{X} \oplus \bold Z H ] = 3$ (In the course of proof of (2.9)); and the fact that $X$ admits no automorphisms of order $7$ if $T_{X}$ is isomorphic to 
$\pmatrix 2 & -1 \\
          -1 & 2
\endpmatrix$ [Vi]. 
Further details are left to the readers. \qed 
\enddemo
Let us next consider the irreducible decomposition of the natural linear 
action of $G$ on $S_{X} \otimes \bold C$. We adopt the notation in [ATLAS]. 
We denote 
the irreducible representation of the character $\chi_{i}$ in [ibid.] 
by $V_{i}$. Then, 
$\text{dim}V_{i}$ is $1$, $3$, $3$, $6$, $7$, $8$ if $i = 1, 2, 3, 4, 5, 6$. 
\proclaim{Claim 3} The irreducible decomposition of the natural action of $G$ on $S_{X} \otimes \bold C$ is: 
$$S_{X} \otimes \bold C = V_{1} \oplus V_{4} \oplus V_{4}' \oplus V_{5},$$ 
where $V_{4}'$ is a copy of $V_{4}$. 
\endproclaim 
\demo{Proof} Set $S_{X} \otimes \bold C = \oplus_{i = 1}^{6} V_{i}^{\oplus n_{i}}$. Here we have $n_{1} = 1$ by (2.1)(1) and also $n_{2} = n_{3}$ as $V_{2}$ 
and $V_{3}$ are (complex) conjugate to each other. Then by counting the dimension, one has 
$$20 = n_{1} + 6n_{2} + 6n_{4} + 7n_{5} + 8n_{6}.$$ 
Recall that for $g \in G$, one has 
\roster 
\item 
$\chi_{\text{top}}(X^{g}) = 2 + \text{tr}(g^{*} \vert S_{X}) + \text{tr}(g^{*} \vert T_{X}) = 4 + \text{tr}(g^{*} \vert S_{X})$ 
\item $\chi_{\text{top}}(X^{g})$ is $8$, $6$, $4$, $3$ if $\text{ord}(g) = 
2$, $3$, $4$, $7$. 
\endroster
The first equality is nothing but the Lefschetz fixed point formula and the second one is the result of Nikulin [Ni]. By applying these two formula and the character table in [ATLAS] for elements of $G$ of order $2$, $3$, $4$, $7$ respectively, one obtains:
$$n_{1} - 2n_{2} + 2n_{4} - n_{5} = 4,$$
$$n_{1} + n_{5} - n_{6} = 2,$$
$$n_{1} + 2n_{2} - n_{5} = 0,$$
$$n_{1} - n_{2} - n_{4} + n_{6} = -1.$$ 
Combining all of these together, we find that $n_{1} = 1$, $n_{2} = n_{3} = 0$, $n_{4} = 2$, $n_{5} = 1$, $n_{6} = 0$. This gives the result. \qed 
\enddemo 
Let $\tau$ be a generater of $\mu_{3}$. We regard $\tau \in \tilde{G}$ through 
the isomorphism found in Claim 1. Since $\tau^{-1}g\tau = g$ for all $g \in G$ 
and $\tau$ is of order $3$, one has $\tau(V_{i}) = V_{i}$ for each $i$ and 
$\tau(V_{4}') = V_{4}'$. Then by Schur's Lemma, $\tau \vert V_{i}$, $\tau \vert V_{i}'$ are all scalar multiplications. Note that $\tau \vert V_{1} = 
id.$ for $S_{X}^{\tilde{G}} \not= \{0\}$. Set 
$\tau \vert V_{4} = \zeta_{3}^{a}$, 
$\tau \vert V_{4}' = \zeta_{3}^{b}$, and 
$\tau \vert V_{5} = \zeta_{3}^{c}$, 
where $a, b, c \in \bold Z/ 3$. Since $\tau \vert S_{X} \otimes \bold C$ is 
defined over 
$S_{X}$, the multiplicities of eigenvalues $\zeta_{3}$ and $\zeta_{3}^{2}$ of $\tau \vert S_{X} \otimes \bold C$ are the same. Therefore, $c = 0$ and $a + b = 0$, i.e. $(a, b, c)$ is either $(0, 0, 0)$ or $(1, 2, 0)$. 
\par
Consider first the case $(a, b, c) = (0, 0, 0)$. Then $\tau^{*} \vert S_{X} = id.$ (and $\tau^{*}\omega_{X} = \zeta_{3}\omega_{X}$). However $T_{X}$ would then be a $3$-elementary lattice by [OZ2, Lemma (1.3)], a contradiction to Claim 2.
\par
Let us consider next the case $(a, b, c) = (1, 2, 0)$. Let $g$ be an 
order $2$ element of $G$. Set $h := \tau g \in \tilde{G}$. 
Then $h$ is of order $6$ and satisfies $h^{3} = g$. In particular, 
one has $X^{h} \subset X^{g}$. Here $X^{g}$ is an $8$-point set by [Ni]. 
Thus, $X^{h}$ also consists of finitely many points (possibly empty), 
thereby $\chi_{\text{top}}(X^{h}) \geq 0$. On the other hand, 
by using the Lefschetz fixed point formula, the fact $(a, b, c) = (1, 2, 0)$, Claim 3, and the character table [ATLAS], one calculates
$$\chi_{\text{top}}(X^{h}) = 2 +  \text{tr}(h^{*} \vert S_{X}) + \text{tr}(h^{*} \vert T_{X})$$
$$ = 2 +  \{1 + 
\text{tr}(g \vert V_{4})(\zeta_{3} + \zeta_{3}^{2}) + \text{tr}(g \vert V_{5})\cdot 1 \} + (\zeta_{3} + \zeta_{3}^{2})$$ 
$$ = 2 + 1 + 2\cdot (-1) + (-1)\cdot 1 + (-1) = -1 < 0,$$ 
 a contradiction to the previous inequality $\chi_{\text{top}}(X^{h}) \geq 0$. Now we are done. \qed 
\head Proof of Proposition 2 \endhead 
Let $\Lambda$ be the K3 lattice, i.e. the lattice $U^{\oplus 3} \oplus E_{8}(-1)^{\oplus 2}$. Choosing a marking $\tau : H^{2}(X_{168}, \bold Z) \simeq \Lambda$,  we set $\Lambda_{0} := \tau(H^{2}(X_{168}, \bold Z)^{G})$. This $\Lambda_{0}$ is a positive definite even lattice of rank $3$ whose $\bold R$ linear extension is spanned by the image of the classes of the invariant ample class $\eta$, $\text{Re}(\omega_{X_{168}})$ and $\text{Im}(\overline{\omega}_{X_{168}})$. 
Fixing the Ricci flat K\"ahler metric $g$ on $X_{168}$ such that the cohomology class of the associated $(1,1)$-form is $\eta$ and regarding $\Lambda_{0} \otimes \bold R$ as a HK $3$-space, one obtains the twister family $f : \Cal X \rightarrow \bold P^{1}$ with $\Cal X_{0} = X_{168}$ (See for instance [Be, Expos\'e X]). This $f$ is a smooth 
non-isotrivial family of (not necessarily algebraic) K3 surfaces $\Cal X_{t}$. 
Denote by $\omega_{t}$ a nowhere vanishing holomorphic two form on $\Cal X_{t}$ and by $\eta_{t}$ the K\"ahler class on $\Cal X_{t}$ associated with $g$. Let us fix a marking $\tilde{\tau} : R^{2}f_{*} \bold Z \simeq \Lambda$ such that $\tilde{\tau}_{0} = \tau$. (Here we used the fact that $\bold P^{1}$ is simply-connected.) By the construction, for each $t \in \bold P^{1}$, the HK $3$-space $\Lambda_{0} \otimes \bold R$ is spanned by the three vectors $\tilde{\tau}_{t}(\eta_{t})$, $\tilde{\tau}_{t}(\omega_{t})$ and $\tilde{\tau}_{t}(\overline{\omega}_{t})$. In particular, we have $\rho(\Cal X_{t}) \geq 19$ for all $t \in \bold P^{1}$. There are then infinitely many $t$ such that $\rho(\Cal X_{t}) = 20$ by [Og2]. Such $\Cal X_{t}$ is necessarily algebraic by [SI] and $\tilde{\tau}_{t}(T_{\Cal X_{t}})$ is a primitive sublattice of rank $2$ of $\Lambda_{0}$. Using the marking $\tilde{\tau}$, let us define the (real) period map:
$$\iota \circ p : \bold P^{1} \rightarrow \{[\omega] \in \bold P(\Lambda \otimes \bold C) \vert (\omega, \omega) = 0, (\omega, \overline{\omega}) > 0 \}$$ 
$$\simeq \{T \in \text{Gr}^{+}(2, \Lambda \otimes \bold R) \vert T \, \text{is positive definite}\}.$$
\par 
Since $p$ is a complex analytic map and $\bold P^{1}$ is compact, $\iota \circ p$ is finite. Therefore, for each rank two sublattice $T$ (of $\Lambda_{0}$), there are at most finitely many $t \in \bold P^{1}$ 
such that $\tilde{\tau}_{t}(T_{\Cal X_{t}}) = T$. Hence, by the global Torelli Theorem for K3 surfaces with the maximum Picard number 20 ([SI]), the family $f$ contains infinitely many non-isomorphic algebraic K3 surfaces. Now the following Claim completes the proof of Proposition 2: 
\proclaim{Claim} $\Cal X_{t}$ satisfies $G \subset \text{Aut}(\Cal X_{t})$ 
for all $t \in \bold P^{1}$. 
\endproclaim
This Claim also shows that there are uncountably many (non-algebraic) K3 surfaces admitting $L_{2}(7)$-actions.  
\demo{Proof} Since $G \vert (\tilde{\tau}_{0})^{-1}(\Lambda_{0}) = \{id. \}$ and $\eta_{t}, \text{Re}(\omega_{t}), \text{Im}(\omega_{t}) \in (\tilde{\tau}_{t})^{-1}(\Lambda_{0} \otimes \bold R)$, we see that $(\tilde{\tau}_{t})^{-1} \circ (\tilde{\tau}_{0}) \circ G \circ (\tilde{\tau}_{0})^{-1} \circ (\tilde{\tau}_{t})$ is an effective Hodge isometry of $H^{2}(\Cal X_{t}, \bold Z)$. This action is also faithful, beacuse $G$ is simple and $G \vert (\tilde{\tau}_{0})^{-1}(\Lambda) \not= \{id. \}$. Hence $G \subset \text{Aut}(\Cal X_{t})$ for each $t \in \bold P^{1}$ by the global Torelli Theorem for K3 surfaces. \qed \enddemo  
\head References added \endhead 
\par \noindent
[Be] A. Beauville, Geom\'etrie des surfaces K3, Ast\'erisque 126 (1985). 
\par \noindent
[Og2] K. Oguiso, Picard numbers in a fmily of hyperk\"ahler manifolds - a supplement to the article of R. Borcherds, L. Katzarkov, T. Pantev, N. I. Shepherd-Barron, preprint (2000) submitted. 
\par \noindent
[OZ2] K. Oguiso and D. -Q. Zhang, On Vorontsov's Theorem on $K3$ surfaces 
with non-symplectic group actions, Proc. AMS {\bf 128} (2000)
1571--1580.
\par \noindent
[Ya] S. T. Yau, On the Ricci curvature of a compact K\"ahler manifold and the complex Monge-Amp\'ere equations I, Comm. Pure Appl. Math. 31 (1978) 339 --411.

\enddocument